
\input amstex
\documentstyle{amsppt}
%
%
%
\magnification=1200
\hoffset=0.4cm\voffset=0.4cm
\hsize=13.5cm\vsize=18.5cm
\baselineskip=13pt plus 0.2pt
\parskip=4pt
%
%
%
\def\R{\Bbb R }
\def\N{\Bbb N }
\def\P{\Bbb P }

\def\Z{\Bbb Z }
\def\C{\Bbb C }

%
%
\font\bigbf=cmbx10 scaled 1200

\font\small=cmr8
\TagsOnRight
\NoRunningHeads
\define\np{\vfil\eject}
\define\nl{\hfil\newline}

\redefine\l{\lambda}

\define\a{\alpha}

\redefine\t{\tau}
\redefine\i{{\,\bold i\,}}
\define\PL{Phys\. Lett\. B}
\define\NP{Nucl\. Phys\. B}
\define\LMP{Lett\. Math\. Phys\. }
\redefine\CMP{Commun\. Math\. Phys\. }

\define\FA{Funktional Anal. i. Prilozhen.}
\define\mapleft#1{\smash{\mathop{\longleftarrow}\limits^{#1}}}
\define\mapright#1{\smash{\mathop{\longrightarrow}\limits^{#1}}}
\define\mapdown#1{\Big\downarrow\rlap{
   $\vcenter{\hbox{$\scriptstyle#1$}}$}}
\define\mapup#1{\Big\uparrow\rlap{
   $\vcenter{\hbox{$\scriptstyle#1$}}$}}
\define\cint #1{\frac 1{2\i\pi}\oint_{C^{#1}}}

\define\im{\text{Im\kern1.0pt }}

\define\RS{Riemann surface}
\define\RSs{Riemann surfaces}

\redefine\deg{\operatornamewithlimits{deg}}
\define\ord{\operatorname{ord}}
\define\fpz{\frac {d }{dz}}
\define\fpt{\frac {d }{dt}}

\define\pfz#1{\frac {d#1}{dz}}
\define\pfx#1{\frac {d#1}{dX}}
\define\pft#1{\frac {d#1}{dt}}
\define\KNN {Krichever - Novikov }

%
%
\hfill Mannheimer Manuskripte 137

\hfill January 1992

\vskip 3cm
\centerline{\bigbf Degenerations of Generalized}
\centerline{\bigbf \KNN Algebras on Tori}
\vskip 2cm
\centerline{Martin Schlichenmaier}
\vskip 0.5cm
\centerline{\small Dept. of Mathematics and Computer Science}
\centerline{\small University of Mannheim, A5, P.O. Box 103462}
\centerline{\small D-W-6800 Mannheim 1, Germany}
\centerline{\small e-mail: FM91\@DMARUM8.bitnet}
\vskip 1cm\noindent
{\bf Abstract.}
Degenerations of Lie algebras of meromorphic vector fields
on elliptic curves (i.e\. complex tori) which are holomorphic
outside a certain set of points (markings) are studied.
By an algebraic geometric degeneration process certain subalgebras of
Lie algebras of meromorphic vector fields on $\P^1$ the Riemann sphere
are obtained. In case of some natural choices of the markings
these subalgebras are explicitly determined. It is shown that
the number of markings can change.
\np\pageno=1
{\bf 1. Introduction}
\vskip 0.5cm
The Virasoro algebra is of fundamental interest in
certain branches of mathematics
and  physics (conformal field theory etc.). The Virasoro algebra
is the graded Lie algebra given as the universal central
extension of the algebra of meromorphic vector fields on the
Riemann sphere $S^2$ holomorphic outside the points $z=0$ and
$z=\infty$. The grading is given by the order of the zeros
at the point $z=0$ (a negative order corresponds to a pole).
For details, see Section 5.

In \cite{11} Krichever and Novikov introduced a generalization
of this graded algebra to higher genus \RSs. The algebra (without
central extension) consists of vector fields meromorphic on the
compact Riemann surface $X$ and holomorphic outside two
arbitrary generic
but fixed points. Of course, this step of the generalization is
quite forward. Their important step was to introduce a
grading by giving a
suitable basis of the algebra by defining them to be the
homogeneous elements.
With respect to this grading the Lie algebra is almost graded (see Sec. 3).

In the interpretation of $X$
together with the two points
as the string world sheet of one incoming
and one outgoing string it is quite natural to ask for generalizations
of these almost graded Lie algebras for the case of several strings
coming in and several strings (not necessarily the same number) going out.
Such a generalization was given by me \cite{16-18} (see also
\cite{14,15}). For related work of other authors on such
generalizations see \cite{6,13,9}.
Again the first step is quite forward. The algebra without central
extension consists of vector fields which are meromorphic
on $X$ outside a fixed finite set of points.
After dividing up this set in two nonempty subsets which I call
in- resp. out-points a grading can be introduced in such a way
that the Lie algebra structure is almost graded.
In fact, all of this can be done for the whole Lie algebra of
differential operators of degree less or equal   1, which is a
semidirect product of the
vector field algebra with the abelian algebra of functions.
Generalizations of affine Kac-Moody algebras
can be introduced (see \cite{17,18}).

If one deforms the complex structure of the \RS\ together
with the set of points a lot of interesting problems appear.
For example, it is quite natural to study such  deformations
in string theory.
Here one has to `integrate' over the whole moduli space of
marked \RSs.
This moduli space is not compact.
To compactify it one has to introduce
`singular Riemann surfaces' and to study the behaviour of the
nonsingular objects when approaching the boundary of the moduli space.
After `resolving the singularity' (see below), one obtains
honest \RSs\ of lower genera. By this one obtains a relation
between disjoint unions of moduli spaces for
nonsingular \RSs. (see \cite{8} for possible physical interpretation).

To perform such a study one has to know how the basis of the
vector field algebra or more generally the basis of the
moduls of forms of
arbitrary weights behave under the deformation of the
complex structure of the \RS\ and under the deformation of the
positions of the points where poles are allowed.
For this it is quite useful to have explicit expressions
of the basis elements. Such were given in \cite{15,17}
in terms of rational functions (for genus $g=0$) , theta functions
$(g\ge 1)$ and Weierstra\ss\  $\sigma$ function ($g=1$).
Even in these representations it is very hard to see what happens
under degeneration to singular `Riemann surfaces'.

By using the Weierstra\ss\ $\wp$ function
Deck \cite{3} and Ruffing \cite{12} (see also \cite{5} )
studied the case of the torus (i.e. $g=1$) with 2, 3 and 4 `punctures'
in certain symmetric positions.
They were able to express the structure constants of the vector
field algebra in purely algebraic deformation parameters,
and observed that
for certain values corresponding to certain degenerations  there
is a (formal) relation to the
Virasoro algebra (see also \cite{4}).

In this letter I want to study the deformation for the $g=1$ case in
detail and give a precise formulation on geometrical grounds
of the relation between the
involved objects.
Having at least a large part of the prospective readership in mind, I
have to say a few words on the algebraic geometric language before
I can give a rough outline of the result. This will be done in
Section 2.
\nl
In Section 3 I translate the setting of marked complex tori into
the setting of marked elliptic curves (i.e. nonsingular cubic curves).\nl
In Section 4 I study the degenerations of the marked cubic curves with
some  natural choices of the markings into marked singular
cubic curves, together with the basis of the
associated vector field algebras.\nl
In Section 5 the objects obtained by degenerations are identified
with objects living on the complex projective line (i.e. the
Riemann sphere).
\nl
The main result can be found in Theorem 1 in this section.

In this letter I will restrict myself to the case of two and three
markings. All possible effects occur already  here.
The studies have been done for forms of arbitrary integer weight
$\l$. Due to lack of space I will concentrate on the vector
fields (i.e. $\l=-1$). The reader will have no difficulties doing
the general case including the differential operator algebra
along the outline given in this letter by himself.
\np
{\bf 2. A few words on the language of
algebraic geometry and an outline of the results}
\vskip 0.3cm
The right theory to describe degenerations of curves is
the theory of algebraic geo\-metry. Here we can confine ourselves with the
language of algebraic varieties.
Roughly speaking, a closed variety is
the set of zeros of a polynomial. A (complex) affine or projective
curve is an closed algebraic variety over $\C$ of dimension 1
either in affine or projective space.
I will call them just curves.
For the notions above and in the following  see  \cite{19} for a
quick introduction or any other book on algebraic geometry.
Singular curves are
incorporated in the theory of algebraic curves from the beginning.
By well-known theorems due to Chow and Kodaira every compact
\RS\ $X$ can be holomorphically  embedded as a complex
nonsingular projective curve in projective space in a very strong
sense. For example, meromorphic functions
(differentials)  on $X$ become rational
functions (differentials) on the embedded curve \cite{19,p.60}.
Recall, a rational function on the embedded curve
in $\P^n$ can be given as
quotient of homogeneous polynomials of the same
degree in $n+1$ variables.

Using the differential equation of the associated
Weierstra\ss\ $\wp$ function I will give in Section 3
the well-known embedding for a torus into
the complex plane $\P^2$ in an explicit manner

By this embedding all tori can be identified as cubic curves
(i.e. curves on which the points  are the zeros of
a homogeneous cubic polynomial)
which are nonsingular. Such curves are termed elliptic curves.
By allowing the cubic polynomial to degenerate one obtains two
nonisomorphic types of irreducible singular cubic curves.
The first one is the nodal cubic $\ E_N\ $ with one singular point
which is a node, i.e. a point where two branches of the curve meet
transversly. The other one is the cuspidal cubic $\ E_C\ $ with one
singularity, which is a cusp, i.e. a point with only one branch of
the curve meeting this point with higher multiplicity.
In Section 4 I will give exact polynomials describing  these curves
(see also \cite{19,p.79}).

By an algebraic geometric process (normalization or blow-up)
every singular curve $X$ can be desingularized, i.e. there is a
nonsingular curve $\widehat{ X}$ and a surjective algebraic map
$\ \psi:\widehat{ X}\to X\ $ such that outside the singularities of $X$
$\psi$ is an isomorphism and above the singularities there are only
finitely many points of $\widehat{X}$.
Moreover, the desingularization curve is  up to algebraic
isomorphisms uniquely defined by the above features.

A family of curves consists of 2 (not necessarily closed) varieties
$B$ and $\Cal X$ and a surjective algebraic map $\phi:{\Cal X}\to B$
such that the fibres $\phi^{-1}(b)$ for every base point
$b\in B$ is an algebraic curve. The different fibres are called
the members of the family. $B$ is called the base.
Intuitively, a family of curves consists of curves depending on
algebraic parameters. If $X$ and $X'$ are different curves, but there
is a family $\Cal X$ over an irreducible base
which contains both curves as members,
we call $X$ a deformation of $X'$ and vice versa.

If $X$ is a nonsingular member of  a family
of curves over an irreducible base
then all other nonsingular members have the same genus $g$ as $X$.
For members of the family which are singular curves the genus $g$
 defined as dimension of the
global holomorphic differentials does not make sense.
\footnote{The correct invariant object for the family is the
arithmetic genus $\ p_a=\dim H^1(X,\Cal O_X)\ $ of the curves $X$.
Due to Serre duality it coincides for nonsingular curves
with the (geometric) genus $g$.}
But at least we know that the genus of their desingularization
will be less then  $g$.
Hence, in our situation the desingularization
of $E_N$ and $E_C$  is a nonsingular
curve of genus 0. This only could be $\P^1$ the complex projective
line. We obtain maps
$$\psi_N:\P^1\to E_N,\qquad\psi_C:\P^1\to E_C\ .\tag 2-1$$
For the nodal cubic there are 2 points lying above the singularity,
for the cuspidal cubic there is one point lying above the singularity.

If objects are defined in a global way
for a family of cubic curves (maybe containing singular ones)
  they are defined
for every member (curve) including the singular ones.
In the case I will consider in this letter the meromorphic vector fields
etc. can be expressed as rational expressions in the affine coordinate
functions $X$ and $Y$
of the complex plane.
Considered as functions on  the curves
they will depend
on the algebraic deformation parameters. The
objects will make sense also in the limit on the singular cubics.
The $X$ and $Y$ coordinate functions of the curve in the limit
can be expressed as functions of the affine coordinate
function $t$ on $\P^1$,
i.e. they can be pullbacked to $\P^1$. By this process they define
objects on $\P^1$.

We have to keep in mind that in our family of cubic curves
every curve posesses a certain finite set of marked points
(sometimes called `punctures') where our objects are allowed to
have poles. If we want to speak of a deformation
of marked curves  we have to
require that the marked points will vary `smoothly' along
the different members. In more precise terms, the maps assigning to
every base point of the family the corresponding markings are required to
be algebraic maps.

The most natural choice of points on elliptic curves are the $\ n-$torsion
points ($n\in\N$).
They are defined as follows.
 An elliptic curve $E$ comes with an abelian
group structure compatible with the algebraic geometric structure.
A $n-$torsion point   is a point $a\in E$ with $n\cdot a=0$.
It is called a primitive $n-$torsion point if $n$ is the
smallest such natural number.
In the complex analytic picture tori are realized as
$$T\ =\ \C/L, \qquad L=\langle \ 1\ ,\ \tau\ \rangle, \qquad \tau\in\C,\
\im \tau>0\tag 2-2$$
where $L$ is a lattice. Here the  above group structure is nothing
else as
addition in $\C$. \nl For example the primitive $2-$torsion points are the
point
($\mod L$)
$$w_1=\frac 12,\qquad
w_2=\frac {1+\t}2,\qquad w_3=\frac {\t}2\ .\tag 2-3$$
Such $n-$torsion points can universally be defined for the whole
family of elliptic curves.
Only such markings will be considered here.

Now I come to an outline of the results.
In Section 3 I will define a family of cubic curves depending on
3 Parameters $e_1,e_2$ and $e_3$. The curves are nonsingular if
all 3 parameters are pairwise distinct. By giving different
choices of the marking I obtain different families of
marked cubic curves.
We can continue the whole situation to all possible
values of the parameters $e_1,e_2$
and $e_3$. We obtain above the additional points
singular cubic curves.
In the case that the marking maps  stay distinct
also in the limit case and none of marked
points in the limit case  coincide with a nodal singularity
(cuspidal singularities are allowed) we obtain from an torus
situation with $N$ markings (punctures) a Riemann sphere situation
with the same number of markings.
If one of the marked points coincide with a nodal singularity we
obtain from a $N$ marking situation on the torus an $N+1$ situation
on the sphere.
Of course, it is also possible that the some markings coincide in
the limit. This will give us a transition to situation with less markings.
We see that if we want to consider degenerations we can not stay
at a fixed
number of markings.

Because we have to study forms of weight $\l$, (i.e. sections in the
$\l$-tensor power of the canonical bundle $K_E$) we have to examine the
behaviour of $\ (dz)^\l\ $ under the algebraization and the pullback.
Here an additional effect at the singularities occur.
In the case
of the vector fields ($\l=-1$) additional zeros are introduced.
\footnote {In the case of  the differentials ($\l=1$)  additional
poles of certain orders are introduced.}
Hence we obtain only a subalgebra in the degeneration,
if the singularity is not a marked point.

By this well defined process of algebraic geometric degenerations
we have the  deformations of the algebras of vector fields
on the elliptic curves (resp\. tori) with  $N$ markings
under explicit control.
In the 3 point case  I studied the following degenerations
were obtained:
\nl
(1) the full $3$ point vector field algebra on $\P^1$, or
(2) some
subalgebra which  can explicitly given,
(3) the full Virasoro-algebra (without central extension) or
(4) a explicitly given subalgebra of it.
\nl
Which case occurs
depends on the markings and
the algebraic geometric deformation under consideration.
Theorem 1 in Section 5 gives the exact result.
Such a list can easily given for every $N$.

Let me mention here that it is also possible to study
deformations of the situation  in the sense that the
genus (i.e. the topology) does not change but two marked points
come together (see \cite{5}).
\np
{
\bf 3. From marked tori to marked elliptic curves}
\vskip 0.5cm
By definition a torus $T$ is the quotient of $\C$ by a lattice
$L=\langle \omega_1,\omega_2\rangle\ $, $T=\C/L$, where
$\omega_1,\omega_2\in\C$ are linearly independent over $\R$.
The complex structure of the quotient is the complex structure induced
from $\C$. Up to analytic isomorphy $T$ can be given as quotient with
respect to a lattice $L$ as given in (2-2).
The field of meromorphic functions on $T$ can be
identified with the field of doubly-periodic functions on $\C$, which
in turn can be generated by the Weierstra\ss\ $\wp$ function
and its derivative $\wp'$ subject to the relation
$$\wp'(z,\tau)^2=4\wp(z,\tau)^3-g_2(\tau)\wp(z,\tau)-g_3(\tau)\ .\tag 3-1$$
Here everything depends on the lattice parameter $\tau $.
The discriminant\nl
$\ \Delta(\t)\ =\ g_2(\tau)^3-27g_3(\tau)^2\ $
is for every admissible $\tau$ different from zero.
For further reference we collect the following well-known results.
$\wp$ has a pole of order 2 at the lattice points and is an even function.
$\wp'$ has a pole of order 3 at the lattice points and is an odd
function. It has zeros at the
the primitive $2-$torsion points on $T$.

Introducing $\ e_i(\tau)=\wp(z_i,\tau)\ $ we can reformulate (3-1)
$$\wp'(z,\t)^2=4\big(\wp(z,\t)-e_1(\t)\big)\big(\wp(z,\t)-e_2(\t)\big)
\big(\wp(z,\t)-e_3(\t)\big)\ .\tag 3-2$$
The condition $\ \Delta(\t)\ne 0\ $ is equivalent to the fact,
that alle $e_i(\t)$ are
distinct.
\footnote{for more information see \cite{10,19},etc.}
 We  obtain  $\ e_1(\t)+e_2(\t)+e_3(\t)=0\ $.

Now we want to rewrite  the situation
in terms of algebraic geometry. We embed the torus $T$
given by an arbitrary but fixed $\t$ into
$\P^2$ the complex projective plane via
$$\Phi:T\to\P^2,\quad
z\mod
L \mapsto \cases (\wp(z):\wp'(z):1),&z\notin L\\
(0:1:0)&z\in L\ .\endcases\tag 3-3$$
If we denote the homogeneous coordinates in $\ \P^2\ $ by
$\ (x:y:z)\ $ the image of $\Phi$ is the set of zeros of the
homogeneous polynomial
$$Y^2Z-4(X-e_1(\t)Z)(X-e_2(\t)Z)(X-e_3(\t)Z)\ .
\tag 3-4$$
Hence, it is  an algebraic curve.
which turns out to be nonsingular.
This curve $E$ is called an elliptic curve.
 Indeed $\Phi$ is an analytic isomorphism
on its image. Conversely every zero set of a homogeneous cubic polynomial
which is a nonsingular curve can be  described (after a suitable
coordinate transformation) as a zero set of a polynomial
of the type (3-4) with $e_1,e_2,e_3\in\C$, $e_1+e_2+e_3=0$
and comes from an analytic situation.
Notice, in (3-4) all $e_i$ are on equal footing. This represents the
fact that in the algebraic geometric picture all
primitive 2-torsion points are equivalent.
If we define
$$B:=\{(e_1,e_2,e_3)\in\C^3\mid
e_1+e_2+e_3=0, \ e_i\ne e_j \ \text{for}\ i\ne j\},$$
then $B$ is a open subvariety of the affine (even linear) variety
$$\widehat{B}:=
\{(e_1,e_2,e_3)\in\C^3\mid
e_1+e_2+e_3=0,\  \}\ .$$
Eq. (3-4) defines a family $\Cal X$ of nonsingular cubic curves
over the base $B$. Obviously (3-4) also makes sense over
the whole of $\widehat{B}$. The additional members of the
family are the singular cubic curves which will be dealt
with in Section 4.

In the following the point $0 \mod L$ will always be a marking
and has nothing to do with the singularity. Hence, it is possible and
convenient to use affine coordinates in $\P^2$ (obtained by setting
the 3. coordinate equal to 1). If we denote $(0:1:0)$ by the symbol
$\infty$ we can rewrite (3-3)
$$z\mod L \ \mapsto\  (\wp(z),\wp'(z)),\quad
 z\notin L,\qquad \text{and}\quad
L\ \mapsto\ \infty\ .\tag 3-5$$
The affine part of the elliptic curve is given as the zero set
of the polynomial
$$Y^2-4(X-e_1)(X-e_2)(X-e_3)\ .\tag 3-6$$
Under this identification the affine coordinate
function $X$ corresponds to
$\wp$ and $Y$ corresponds to $\wp'$. The field of meromorphic  function on
the torus corresponds to the field
of rational functions on the curve $E$
$$\C(X)[Y]/\big(Y^2-f(X)\big) ,
\quad f(T)=4(T-e_1)(T-e_2)(T-e_3)\ .\tag 3-7$$
Every meromorphic  function on the torus can be described as
a rational function (i.e. as a function which is a
quotient of polynomials) in $X$ and $Y$. In fact, (3-7) shows even that
it can be given as a rational function in $X$ plus $Y$ times
another rational function in $X$.
A function $f$ is called regular on a subset of the curve $E$
if $f$ can be given as quotient of polynomials
such that the  the denominator
polynomial does not vanish on $E$. Regular functions correspond to
holomorphic functions. Due to this equivalence I will usually
call rational (regular) functions meromorphic (holomorphic) functions.

We have to introduce markings, i.e. choose points where poles are
allowed. As explained in Section 2 the natural choices are
$n-$torsion points.
I will consider here  the following cases.

(Case A). In the  two point case I choose $\ z_0=0\ $ and
$\ z_1=1/2\ $ in the
analytic picture. In the elliptic curve picture this
corresponds to choosing the point $\ \infty\ $
and the point with the affine coordinates $\ (e_1,0)\ $.
A different choice  of a $2-$torsion point instead of $z_1=1/2$
(as done in \cite{2})
yields an isomorphic algebra.

(Case B). In the three point case I choose
$\ z_0=0\ $, $\ z_1=1/2+q\ $, and \nl $\ z_2=1/2-q\ $
with the two subcases
$\ q=\tau/4\ $ and $\ q=1/2+\tau/4\ $.
In both cases $z_1$ and $z_2$ are primitive $4-$torsion points. In fact,
these are the two pairs of solutions to\nl
$\ 2\cdot(z \mod L)=\tau/2\mod L\ $.
Other choices of $4-$torsion points
will give the same behaviour of the above
situations (studied at different deformations).
On the elliptic curve, the markings correpond to the points
$\ \infty\ $, $\ (a,b)\ $, and $\ (a,-b)\ $ with
$$a=\wp(1/2+q),\quad
b=\wp'(1/2+q)=\sqrt{4(a-e_1)(a-e_2)(a-e_3)}\ .\tag 3-8$$
Both $(a,b)$ and $(a,-b)$ occur as markings. Hence,
there is no ambiguity of
sign in the complex square root of $b$ .
We are interested in varying the algebraic parameters
$e_1,e_2$ and $e_3$. For this goal we have to express $a$
(and automatically $b$)
in terms of them. Using the addition theorem of the $\wp$
function (see \cite{7}) we obtain
$$a(e_1,e_2,e_3)\ =\ e_3+\sqrt{e_3-e_1}\sqrt{e_3-e_2}\ .\tag 3-9$$
Here the sign of the product  of the square root
has to be choosen following  the rules given in \cite{7}.
The two possibilities
occuring
correspond to the two choices for $q$.
We obtain two pairs of values   $\ (a_1,b_1)\ $, $\ (a_1,-b_1)\ $ and
 $\ (a_2,b_2)\ $, $\ (a_2,-b_2)\ $.

The main reason for choosing the above markings is that there coincide
with the cases considered by Deck and Ruffing \cite{3-5,12}.
Hence, I am able to refer to  their
results and avoid redoing the calculation.
By setting $q=0$ in  (Case B) we obtain formally the 2 point case from the
3 point case.
This allows an elegant compact form of notation if we define
$\ (a_0,b_0)=(e_1,0)\ $
as we will see immediately.

\proclaim{Proposition 1}
A basis of the vector space of
meromorphic (rational) functions on the elliptic curve
$E$ holomorphic outside the points $\infty$ and $(e_1,0)$ in the
2 point case ($s=0$), resp. $\infty$, $(a_s,b_s)$, and $(a_s,-b_s)$
for the 3 point cases ($s=1,2$) is given by the elements
$$A_{2k}^s:=(X-a_s)^k,\qquad
A_{2k+1}^s:=\frac 12Y(X-a_s))^{k-1}\ \tag 3-10$$
where $k$ runs over all integers.
($s=0,1,2$ denotes the different cases.)
\endproclaim
\demo{Proof}
Up to reindexing and rewriting it is just the basis introduced by
Deck and Ruffing. For a proof see \cite{5}.
In fact, this is quite easy to see directly. Using the information
on the functions on $E$ following Eq. (3-7) we see that the
denominator polynomial in the rational functions in $X$ could
only be powers of $(X-a_s)$. By developing the numerator polynomial
in powers of $(X-a_s)$ we obtain the result.
\qed \enddemo
Using the following facts
about the orders at the points of the elliptic curve
\nl ($\gamma\in\C,\ a\ne e_1,e_2,e_3$)
$$\gathered
\ord_{\infty}(X-\gamma)=-2,\qquad\ord_{\infty}(Y)=-3,\\
\ord_{(e_i,0)}(X-e_i)=2,\qquad\ord_{(e_i,0)}(Y)=1,\\
\ord_{(a,b)}(X-a)=1,\qquad\ord_{(a,-b)}(X-a)=1\ ,
\endgathered\tag 3-11$$
we obtain
\proclaim{Proposition 2}
The divisors
for the basis elements
are given as follows:

\noindent (a) in the two point case:
$$\align (A_{2k}^0)&=-2k[\infty]+2k[(e_1,0)],\\
(A_{2k+1}^0)&=-(2k+1)[\infty]+(2k-1)[(e_1,0)]+1[(e_2,0)]+1[(e_3,0)],
\endalign$$
(b) in the three point cases: ($s=1,2$)
$$\align(A_{2k}^s)&=-2k[\infty]+k[(a,b)]+k[(a,-b)],\\
(A_{2k+1}^s)&=-(2k+1)[\infty]+(k-1)[(a_s,b_s)]+(k-1)[(a_s,-b_s)]
\\&\qquad+1[(e_1,0)]+1[(e_2,0)]+1[(e_3,0)].\endalign$$
\endproclaim
Recall that a divisor of a function (or generally of a
section of a line bundle) denotes the points where zeros
of this function occur with their muliplicities. A zero
of negative multiplicity is a pole.

In the genus 1 case the canonical bundle $K$ and hence all its
tensor powers are trivial.
If $f$ is a meromorphic function on the torus
$T$ then $f(z)(dz)^{\l}$ is a
meromorphic section of the bundle
$\ K^{\otimes\l}=K^\l$. We have to express the analytic differential $dz$
as
$$dz=\frac {dX}Y\ .\tag 3-12$$
Observe that $(X-\gamma)$ for $\gamma\notin\{e_1,e_2,e_3\}$ is a
uniformizing variable at $(\gamma,\pm\sqrt{f(\gamma)})$
and $Y$ does not vanish there.
For $\gamma=e_1,e_2$ or $e_3$ the function
$(X-\gamma)$ vanishes of second order
at $(e_i,0)$ and hence compensates for the pole of
$\ 1/Y\ $ at this point.
For the tensor powers we obtain
$$(dz)^{\l}=Y^{-\l}(dX)^{\l},\qquad
\fpz=Y\frac{d}{dX}\tag 3-13$$
where we used for the special case $\l=-1$ the usual derivation notation.
Hence a basis of the sections in $K^\l$ with the same regularity
condition as in Prop.~1 are given by the elements
$$A_m^sY^{-\l}(dX)^{\l},\quad m\in\Z\ .\tag 3-14$$
In the case of vector fields we obtain
from Prop.~1
\proclaim{Proposition 3}
A basis of the space of meromorphic (rational)
 vector fields on the elliptic curve
$E$ holomorphic (regular) outside the points $\infty$ and $(e_1,0)$ in the
2 point case (corresponding to $s=0$), resp.
 $\infty$, $(a_s,b_s)$, and $(a_s,-b_s)$
for the two 3 point cases ($s=1,2$) is given by the elements
$$V_{2k}^s:=(X-a_s)^kY\frac{d}{dX},\qquad
V_{2k+1}^s:=\frac 12f(X)(X-a_s))^k\frac{d}{dX}\ \tag 3-15$$
where $k$ runs over all integers  and
$f(X)=4(X-e_1)(X-e_2)(X-e_3)$.
\endproclaim \noindent
By defining
$$\gather\deg(A_n^0Y^{-\l}(dX)^\l)=n,\tag 3-16\\
\deg(A_{2k}^sY^{-\l}(dX)^\l)=
\deg(A_{2k+1}^sY^{-\l}(dX)^\l)=k\tag 3-17\endgather$$
we introduce a grading in the vector space of all sections.

The  vector fields with the above regularity conditions
define a Lie algebra under the usual Lie bracket.
It is usually called generalized Krichever - Novikov algebra
(see \cite{16-18} for details).
The sections in $K^\l$ (again with the same regularity
conditions) are Lie modules over this algebra
with respect to taking the Lie derivative.
The function algebra itself operates as multiplication on the
space of sections. Putting this two actions together
the space of sections are Lie modules over the
Lie algebra of differential operators of degree
$\le 1$, which is given as semidirect product
of the algebra of vector fields with the functions.
Because, the situation is completely analogous in the general
situation I will only do here  the vector field case.
To avoid cumbersome notation I suppress the superscript $s$
\proclaim{Proposition 4}
The Lie algebra structure of the \KNN\ algebra is given by the
following structure  equations in the generators introduced
in Prop. 3, ($k,l\in\Z$)

\noindent (a) in the two point case:
$$\aligned
[V_{2k},V_{2l}]&=(2l-2k)V_{2(k+l)+1}\\
[V_{2k+1},V_{2l+1}]&=((2l+1)-(2k+1))\{
V_{2(l+k+1)+1}+3\,e_1V_{2(l+k)+1}\\&\qquad\qquad
+(e_1-e_2)(e_1-e_3)\,V_{2(l+k-1)+1}\}\\
[V_{2k+1},V_{2l}]&=(2l-(2k+1))
V_{2(l+k+1)}+(2l-(2k+1)+1)\,3\,e_1V_{2(l+k)}\\&\qquad
+(2l-(2k+1)+2)(e_1-e_2)(e_1-e_3)\,V_{2(l+k-1)}\
 ,\endaligned\tag 3-18$$
(b) in the three point case:
(where $a$ denotes $a_1$ or $a_2$ depending on the case $s$ )
$$\aligned
[V_{2k},V_{2l}]&=(2l-2k)V_{2(k+l)+1}\\
[V_{2k+1},V_{2l+1}]&=((2l+1)-(2k+1))\{
V_{2(l+k+1)+1}+3\,a\,V_{2(l+k)+1}\\&\qquad\qquad
+(3a^2-(e_2^2+e_2e_3+e_3^2))\,V_{2(l+k-1)+1}
+(1/4)\,b^2V_{2(l+k-2)+1}\\
[V_{2k+1},V_{2l}]&=(2l-(2k+1))
V_{2(l+k+1)}+(2l-(2k+1)+1)\,3\,a\,V_{2(l+k)}\\&\qquad
+(2l-(2k+1)+2)\big(3a^2-(e_2^2+e_2e_3+e_3^2)\big)V_{2(l+k-1)}
\\&\qquad\qquad+(2l-(2k+1)+3)(1/4)\,b^2\,V_{2(l+k-2)}
 .\endaligned\tag 3-19$$
\endproclaim
\demo{Proof}
This is a reformulation of results
obtained by Ruffing and Deck
which again could be calculated directly by using informations about
the possible pole orders of the Lie bracket.
\cite{5,3,12}.
\enddemo
If $d$ denotes the sum of the degrees of the generators
on the left hand side of (3-18),\nl (3-19) then the degrees of the
generators appearing on the
right hand side are within the range $d-3$ and $d+1$.
Hence, the Lie algebra structure is almost graded with respect to the
gradings (3-16),(3-17).  The same is true for the
Lie module structures of the space of  sections. These
features are important to construct semi-infinite wedge representations
(see \cite{16,17}).

Of course, it would be possible to avoid completely the use of the
complex analytic tori and to start with the cubic curves.
But to see the relation with
the in conformal field theory more familiar analytic picture
I decided not to do so.
To a certain extend an  algebraic formulation
has been given in a recent  preliminary version of a preprint
by Anzaldo-Meneses \cite{1}.
\np
{\bf 4. From marked elliptic curves to marked singular cubic curves}
\vskip 0.5cm
The singular members of the family of cubic curves
 given by the polynomial  (3-4)
are  exactly the curves
lying above points in $\widehat{B}$ where at least two of the $e_i$
 coincide.
We obtain two different situations.\nl If only two coincide
we get the nodal cubic $E_N$ given
by
$$Y^2=4(X-e)^2(X+2e)\ .\tag 4-1$$
Here $e$ denotes the value of the coinciding ones. By $e_1+e_2+e_3=0$
the third one has the value $-2e$.\nl
If all three have the same value, which necessarily equals
$0$, we get the cuspidal cubic $E_C$
$$Y^2=X^3\ .\tag 4-2$$
The singularity on $E_N$ is the point $(e,0)$, a node.
The singularity on $E_C$ is the point $(0,0)$, a cusp.
In both cases the point $\infty=(0:1:0)$ lies on the cubic curves but
is not a singularity.

As explained in Section 2 in both cases the complex projective
line $\P^1$ is the desingularization. The points of $\P^1$ are given
by homogeneous coordinates $(t,s)$. We define the following maps
$\quad \psi_n,\psi_C:\P^1\to\P^2$ by
$$\align
\psi_N(t:s)&=(\ t^2s-2es^3\ :\ 2t(t^2-3es)\ : \ s^3\ ),\tag 4-3\\
\psi_C(t:s)&=(\ t^2s\ :\ 2t^3\ : \ s^3\ )\ .\tag 4-4\endalign$$
Under these maps the point $\ \infty=(1:0)\ $ on $\P^1$ corresponds
to (and only to) the point $\infty$ on both curves $E_N$ and $E_C$.
The maps are given by homogeneous polynomial. Hence, they are
obviously algebraic maps. Again it is enough to consider the
affine part (i.e. we are setting $s=1$). We use the same symbols for
the affine maps.
A direct calculation shows
the following facts.\nl (1) Image $\psi_N=E_N$ and Image
$\psi_C=E_C$. \nl
(2) $\psi_C$ is $1:1$.\nl
(3) The only points where $\psi_N$ is not $1:1$ are the points
$\ t=\sqrt{3e}\ $ and $\ t=-\sqrt{3e}\ $. They both project onto the
singular point $(e,0)$. The point $(-2e,0)$ correspond to
$t=0$.
Hence,
\proclaim{Proposition 5}
The maps
$$\psi_N:\P^1\to E_N\quad\text{and}\quad\psi_C:\P^1\to E_C$$
are the unique desingularization of the nodal cubic $E_N$ resp\.
the cuspidal cubic $E_C$.
\endproclaim
Now we have to consider what happens to the markings.
In all cases $t=\infty$ on $\P^1$ corresponds to a point where
poles are allowed. \nl
Recall that the other markings are given by $\ (e_1,0)\ $ in the
2 point case, resp\. $(a,b)$ and $(a,-b)$ in the 3 point case
where $\ a\ $ and $\ b\ $ are given by (3-8), (3-9).

The situation in the cuspidal case is easy to describe.
Because $e_1=e_2=e_3=0$ all remaining markings go to
the singularity (either 2 point case or 3 point case).
What remains is $\P^1$ with the markings $t=0$ and $t=\infty$.

In the case of nodal degenerations and 2 points
we have to distinguish 2 different degenerations.
\nl (1) If we have $e=e_1$
then the marking becomes the  singular point. Hence on $\P^1$
two points $\ t_1=\sqrt{3e}\ $ and
$\ t_2=-\sqrt{3e}\ $ correspond to the second marking.
We obtain $\P^1$ with 3 markings.
\nl (2) If $e\ne e_1$ then
the marking is an ordinary point. Only the
 point $t=0$ correspond to the second marking.
We obtain $\P^1$ with 2 markings.

In the case of nodal degenerations and 3 points
we have the following possibilities.\nl
(3) If $e=e_1=e_3$ or $e=e_2=e_3$
we find $a_s=e$ and $b_s=0$ for $s=1,2$ using (3-8), (3-9).
Here the two markings $(a_s,b_s),\ (a_s,-b_s) $
join at the singularity. This corresponds to the markings
$t_1=\sqrt{3e}$ and $t_2=-\sqrt{3e}$.
We obtain $\P^1$ with 3 markings.
\nl
(4) If $e_1=e_2$ we have to check the sign of the square root in
(3-9). For $s=1$, (i.e. $q=\t/4$
we obtain $\ a_1=e_1=e\ ,(b_1=0)\ $. Hence the same situation
as in (3). The sign of the product of the roots can be determined
either by following the prescription given in \cite{7}
or by observing that the marking $1/2+\t/4$ is `squeezed'
between the two approaching points $1/2$ and $1/2+\t/2$
under this degeneration.\nl
(5) In the remaining case $e=e_1=e_2$ and
$s=2$ (i.e. $q=1/2+\t/4$) we have to take the other sign of the
product of the roots and obtain $a_2=-5e$ with two associated
values for $b_2$. The markings $(a_2,b_2),\ (a_2,-b_2)$
stay distinct and none of them coincides with the singularity.
They are given by the $t$ values $\pm\i\sqrt{3e}$.
We obtain $\P^1$ with 3 markings.

As shown in Section 3 a section
of $K^{\l}_E$ corresponding to the different situations $s$
for $E$ a
nonsingular cubic curve is given by
\footnote{The ${}' $ indicates that  only finitely many
summands occur.}
$$w(X,Y)=\sum_{k\in\Z}{}'\big(\alpha_k(X-a_s)^k+
\beta_kY(X-a_s)^{k-1}\big)\big(\frac {dX}Y\big)^\l,\qquad
\alpha_k,\beta_k\in\C\ .\tag 4-5$$
On the nonempty open set of nonsingular
points these elements make perfectly sense also in  the degenerate
cases $E_N$ and $E_C$. By pulling it back via  $\psi_N$, resp\. $\psi_C$
on the  open set of $\P^1$ mapping to the nonsingular points we
obtain a well defined meromorphic (rational) Section of
$K^\l_{\P^1}$
$$ \psi_N^*(w)(t)=w(X(t),Y(t))\qquad
\psi_C^*(w)(t)=w(x(t),y(t))\ . $$
To find their representation functions
we have to plug in the expression for $X$ and $Y$ in terms of $\ t\ $ and
to rewrite the differential  in terms of $\ dt\ $.
For the power of the differential we obtain in the nodal case
$$\left(\frac {dX}Y\right)^\l= \frac 1{(t^2-3e)^\l}(dt)^\l,\tag 4-6
$$
and$$
\left(\frac {dX}Y\right)^\l= \frac 1{t^{2\l}}(dt)^\l,\tag 4-7
$$
in the cuspidal case.
By this we see very explicitly  that
additional poles and zeros (depending on the sign of $\l$) are introduced
at the singular points.
A important example is the pullback of the constant
differential $dz=\dfrac {dX}Y$.
We obtain
$$\frac 1{t+\sqrt{3e}}\, \frac 1{t-\sqrt{3e}}\ dt\qquad
\text{resp\.}\qquad
\frac 1{t^2}\, dt\ .$$
We get by degeneration to the nodal cubic a meromorphic
differential with poles of order 1 at the points $\ t=\pm\sqrt{3e}\ $,
 resp\.
a pole of order 2 at $t=0$. These are the Rosenlicht differentials
\cite{20}.

In this letter we are especially interested in the case $\l=-1$, the
vector field case. Here only additional zeros at the the points
above the singularities are introduced.
All Lie derivatives can be calculated on the
(Zariski-) open set
of points which do not become singular points
under the degeneration. Hence, its Lie structure is not inflicted by the
degeneration. We obtain
\proclaim{Proposition 6}
Under the geometric process of degeneration and desingularization
one obtains a geometrically induced `deformation'
of the algebra of vector fields on the torus into a
subalgebra of the vector fields on $\P^1$ consisting of vector
fields with poles only at the points $t$ with $(X(t),Y(t))$ is
a marked point on the degenerate cubic curve.
The degenerate algebra can be obtained by replacing the
generators by their pullbacks and setting the
structure
constants to their limit values.
\endproclaim
By examples,
I will show that the number of markings can change (even increase)  and
that not necessarily the full vector field algebra will occur.

Because it is clear from the description in which situation we are
I drop the index $s$ to avoid cumbersome notations.
Let me give first the pull back of the generators (3-14).
They are ($k\in\Z$)
$$\align
\psi_N^*\big(A_{2k}\left(\frac {dX}{Y}\right)^\l\big)
&\ =\ (t^2-3e)^{-\l}(t^2-2e-a)^k(dt)^\l,\tag 4-8\\
\psi_N^*\big(A_{2k+1}\left(\frac {dX}{Y}\right)^\l\big)
&\ =\ t\,(t^2-3e)^{1-\l}(t^2-2e-1)^{k-1}(dt)^\l,\tag 4-9\\
\psi_C^*\big(A_{2k}\left(\frac {dX}{Y}\right)^\l\big)
&\ =\ t^{-2\l}(t^2-a)^k(dt)^\l,\tag 4-10\\
\psi_C^*\big(A_{2k+1}\left(\frac {dX}{Y}\right)^\l\big)
&\ =\ t^{3-2\l}(t^2-a)^{k-1}(dt)^\l,\tag 4-11
\endalign
$$
The above formula are valid for arbitrary  values of $a$.
Here we only consider the natural choices described above.\nl
In the cuspidal case we see that the limit value of $a$ equals zero.
Hence (4-10), (4-11)  specializes to (in the case $\l=-1$)
$$
\psi_C^*\big(A_{2k}\left(\frac {dX}{Y}\right)^{-1}\big)
\ =\ t^{2k+2}\fpt,\qquad
\psi_C^*\big(A_{2k+1}\left(\frac {dX}{Y}\right)^{-1}\big)
\ =\ t^{2k+3}\fpt\ .\tag 4-12$$
By this we see that in the cuspidal degeneration (with the above
markings) we always obtain the full (2 point) Virasoro algebra.
This coincides with the degeneration of the structure
equations (3-18) and (3-19) to
$$[V_n,V_m]=(m-n)V_{n+m+1}, \quad n,m\in\Z \ .\tag 4-13$$
In the nodal degeneration case, we have to distinguish several
subcases. In the 2 marking case $(a=e_1)$ we obtain the following picture.
\nl (1) For $e=e_1$ the limit marking is a singularity.
For the pullback we obtain
$$\aligned
\psi_N^*\big(A_{2k}\left(\frac {dX}{Y}\right)^{-1}\big)
\ &=\ (t+\sqrt{3e})^{k+1}(t-\sqrt{3e})^{k+1}\fpt,\\
\psi_N^*\big(A_{2k+1}\left(\frac {dX}{Y}\right)^{-1}\big)
\ &=\ t\,(t+\sqrt{3e})^{k+1}(t-\sqrt{3e})^{k+1}\fpt\ .
\endaligned\tag 4-14$$
We see explicitly that on the desingularization
3 markings  occur. From a 2 point situation
we came to a 3 point situation.
The limit algebra is defined by the limit of the structure
equations (3-18) for $e_1=e_2=e$  or $e_1=e_3=e$
$$\aligned
[V_{2k},V_{2l}]&=(2l-2k)V_{2(k+l)+1},\\
[V_{2k+1},V_{2l+1}]&=((2l+1)-(2k+1))\{V_{2(k+l+1)+1}
+3\,e\,V_{2(l+k)+1}\},\\
[V_{2k+1},V_{2l}]&=(2l-(2k+1))V_{2(k+l+1)}
+(2l-(2k+1)+1)3\,e\,V_{2(l+k)}\ .
\endaligned\tag 4-15$$
In the next Section I will identify this
and the following algebras in
more detail.
\nl
(2) The next subcase for the two point situation is
$e_1\ne e$. The limit marking is not a singularity. The pullbacks are
$$\aligned
\psi_N^*\big(A_{2k}\left(\frac {dX}{Y}\right)^{-1}\big)
\ &=\ t^{2k}(t+\sqrt{3e})(t-\sqrt{3e})\fpt,\\
\psi_N^*\big(A_{2k+1}\left(\frac {dX}{Y}\right)^{-1}\big)
\ &=\ t^{2k-1}(t+\sqrt{3e})^2(t-\sqrt{3e})^2\fpt\ .
\endaligned\tag 4-16$$
Only at two points  poles can occur ($t=0$ and $t=\infty$).
Hence, we do not leave the 2 point situation. But
we do not obtain the full Virasoro algebra, because
the vector fields are forced to have zeros at
$t=\pm\sqrt{3e}$.
Indeed, there is a difference in the two formulas of (4-16). We get an
additional condition to identify the exact subalgebra.
This is due to the fact that only those functions on $\P^1$ (given by
Laurent polynomials) can be obtained by pullback from $E_N$ if
they fullfill $\ f(\sqrt{3e})=f(-\sqrt{3e})\ $. For the even functions
this is automatically. For the odd function this forces
them to be zero at the two points. Hence, only functions
generated by $\ t^{2k}\ $ and $\ f^{2k+1}(t^2-3e)\ $ for $k\in\Z$
could occur. (Remember the second term $t^2-3e$ comes from the
differential.)
\nl
The structure equations of the algebra
(obtained by the usual process) are
$$\aligned
[V_{2k},V_{2l}]&=(2l-2k)\,V_{2(k+l)+1},\\
[V_{2k+1},V_{2l+1}]&=((2l+1)-(2k+1))\{V_{2(k+l+1)+1}
-6\,e\,V_{2(l+k)+1}\}+9\,e^2V_{2(l+k-1)+1}\},\\
[V_{2k+1},V_{2l}]&=(2l-(2k+1))V_{2(k+l+1)}
+(2l-(2k+1)+1)(-6\,e)V_{2(l+k)}\\
&\qquad+
(2l-(2k+1)+2)\,9\,e^2V_{2(l+k-1)} .
\endaligned\tag 4-17$$
Now we are coming to the three point cases.\nl
(3)
For $e=e_1=e_3$ or $e=e_2=e_3$
we obtain for both values of $q$
that
the markings $(a,b)$ and $(a,-b)$ move into the
singularity $(e,0)$. This yields exactly the same
generators and the same deformed algebra as
as considered in (1), resp\. in (4-14),(4-15).
Notice, we obtain the Equations (3-18) formally by setting
$a=e_1$ and $b=0$ in (3-19) using the fact $e_1+e_2+e_3=0$.
\nl
(4) For $q=\tau/4$ and the remaining case $e=e_1=e_2$
the situation is the same as described under (3).
\nl
(5) For $q=1/2+\tau/4$ the remaining case is $e=e_1=e_2$.
We had obtained above $a=-5e$ and two different associated
$b$ values.
The pullbacks of the generators are
$$\aligned
\psi_N^*\big(A_{2k}\left(\frac {dX}{Y}\right)^{-1}\big)
\ &=\ (t+\i\sqrt{3e})^k(t-\i\sqrt{3e})^k(t+\sqrt{3e})(t-\sqrt{3e})\fpt\\
\psi_N^*\big(A_{2k+1}\left(\frac {dX}{Y}\right)^{-1}\big)
\  &=\ t\,(t+\i\sqrt{3e})^{k-1}(t-\i\sqrt{3e})^{k-1}
(t+\sqrt{3e})^2(t-\sqrt{3e})^2\fpt\ .
\endaligned\tag 4-18
$$
Poles are at $\ t=\i\sqrt{3e},\ -\i\sqrt{3e},\infty\ $.
 Hence we stay in a three
point situation. As in (2) additional zeros
at $\ t=\pm\sqrt{3e}\ $ are introduced
adn we get an additional condition yielding the difference between
even and odd elements.
We
get only a subalgebra of the corresponding 3 point algebra
on $\P^1$. The structure equations can be easily be obtained from (3-19).
I ommit writing them down here.
\vskip 0.5cm
{\bf 5. The involved algebras on $\P^1$ with arbitrary markings}
\vskip 0.5cm
The Virasoro algebra $\ \Cal V\ $
(without central extension) can be given as the
complex vector space generated by the vector fields
$$L_n\ :=\ t^{n+1}\fpt,\qquad n\in\Z\ .\tag 5-1$$
Its structure equations are
$$[L_n,L_m]=(m-n)L_{n+m},\quad m,n\in\Z\ .\tag 5-2$$
By the geometric reasoning we obtain in the cuspidal cases
for the 2 or 3 point situation with (4-12)
$$\psi_N^*\big(A_n\left(\frac {dX}Y\right)^{-1}\big)=L_{n+1}\ .\tag 5-3$$
Of course, the formal isomorphism (without giving the
geometric explaination) could be seen by the map
$\ \Phi:V_n\mapsto L_{n+1}\ $ in the degenerate structure
equations
(4-13).
\footnote{ Such kind of observations made by Deck and Ruffing
where the starting point of this work. Hence I am very much indebted
for their calculations.}

We are now considering certain subalgebras of the Virasoro algebra.
The vector fields vanishing at the points $\a$ and $-\a$,
with  $\a\ne 0,\infty$
are the vector fields
$$f(t)\,(t^2-\a^2)\,\fpt\tag 5-4$$
where $f(t)$ is a Laurent polynomial in $t$. Obviously, they
define a subalgebra of $\Cal V$ with (vector space) basis
$$t^k\,(t^2-\a^2)\fpt\ = \ L_{k+1}-\a^2L_{k-1},\qquad k\in\Z\ .\tag 5-5$$
Inside this subalgebra we consider the subspace generated by
the vector fields
$$\align
M_{2k}&:=t^{2k}(t^2-\a^2)\fpt=L_{2k+1}-\a^2L_{2k-1},\tag 5-6\\
M_{2k+1}&:=t^{2k-1}(t^2-\a^2)^2\fpt=
L_{2k+2}-2\a^2L_{2k}+\a^4L_{2k-2}\ .\tag 5-7
\endalign$$
Either by direct calculation of the Lie bracket, or calculation
inside the Virasoro algebra we see this is a subalgebra $\ {\Cal W}^\a\ $
of $\ \Cal V$.
\nl For $\a=\sqrt{3e}$ the $M_n$ coincide with the pullbacks (4-16).
The map
$\Phi: V_n\mapsto M_n,\ n\in\Z$ gives the
identification of the geometrically induced degenerated algebra
(4-17) with the algebra ${\Cal W}^\a$.
In \cite{4} a relation of the `deformed' algebra
in the 2 point case  with the Virasoro
algebra was found by purely formal algebraic investigations. Indeed, it
was tried (without sucess) to identify it with the whole
Virasoro algebra. By the above, we see the reason why this could
not work.

I am now coming to the algebra $\ {\Cal Z}^\a\ $  of vector fields on $\P^1$
holomorphic outside the points
$\a,-\a$ and $\infty$ ($\alpha\ne 0,\ \infty$).
In \cite{16,17} a general method for obtaining a basis
was given. If we split the marking into two disjoint sets, the
in-points $\{\a,-\a\}$ and the out-point $\{\infty\}$
the algebra is generated (as vector space)
by elements $e_{n,1}$ and $e_{n,2}$ with
$n\in\Z$. Their divisors are
$$\align
(e_{n,1})\ &=\ (n-1)[\a]+n[-\a]+(-2n+3)[\infty],\\
(e_{n,2})\ &=\ n[\a]+(n-1)[-\a]+(-2n+3)[\infty]\ .\endalign$$
These are the generators of degree $n$. We can symmetrize the
situation by taking suitable linear combinations of them
to obtain generators $H_n$ and $G_n$ with the divisors
$$\aligned
(H_n)&=(n-1)[\a]+(n-1)[-\a]+(-2n+4)[\infty],\\
(G_n)&=(n-1)[\a]+(n-1)[-\a]+1[0]+(-2n+3)[\infty]\ .
\endaligned\tag 5-8$$
These  generators are still homogeneous elements
of degree $n$ and a basis of the algebra ${\Cal Z}^\a$.
They can explicitly be given as
$$
H_n(t):=(t-\a)^{n-1}(t+\a)^{n-1}\fpt,\qquad
G_n(t):=t\,(t-\a)^{n-1}(t+\a)^{n-1}\fpt\ .\tag 5-9$$
A direct calculation shows
$$\aligned
[H_n,H_m]&=2(m-n)G_{n+m-2},\\
[G_n,G_m]&=2(m-n)(G_{n+m-1}+\a^2G_{n+m-2}),\\
[G_n,H_m]&=(2(m-n)-1)H_{n+m-1}+2(m-n)\,\a^2H_{n+m-2}
\ .\endaligned\tag 5-10$$
For $\a=\sqrt{3e}$ the elements (5-9) coincide
with the pullbacks (4-14) obtained in the nodal degeneration in
Section 4 in the cases (1) and (3).
Under the map\nl
$\Phi:\ V_{2k}
\mapsto H_{k+2},\quad
V_{2k+1}
\mapsto G_{k+2}\ $
we obtain the geometrically induced identification of
the degenerate algebra in the case (4-14) with the full
full algebra ${\Cal Z}^\a$.

In the remaining case (5) we have to consider again the subalgebra of
${\Cal Z}^\a$ (now for $\alpha=\i\sqrt{3e}$)
generated by the vector fields with
zeros at
$\pm\sqrt{3e}$ and a similar additional condition as
in the case  ${\Cal W}^\a$.
Again this algebra can be identified with the degenerated algebra
obtained in Section 4. The results are completely
analogous
hence I will not write them down here.

We should not forget that it is not only the algebras, but also their
grading which is of importance. The only graded Lie algebra (in the
strict sense) here is the full Virasoro algebra.
All other algebras are only almost graded.
Hence it does not make sense to investigate whether the maps $\Phi$
are grading preserving. Instead I introduce the following

\noindent{Definition.}
{\sl
Let $T$ and $S$ be almost graded Lie algebras. A Lie homomorphism
\nl $\Phi:S\to T\ $  respects the almost grading
if there are positive natural numbers $\ a,k,l\ $ such that
for every homogeneous element $s\in S$ we get
$$a\cdot\deg (s)-k\ \le\ \deg(\Phi(s))\ \le\  a\cdot\deg (s)+l\ .$$
}

\noindent
By inspection of the formalas above and
collecting the results we obtain
\proclaim{Theorem 1}
Under the geometrically defined deformation of the 2 or 3 point
vector field algebra on the ellitic curve (i.e. torus) with
the markings $\infty$ and $(e_1,0)$ in the 2 point case, and
$\infty$, $(a_s,b_s)$ and  $(a_s,-b_s)$ for the two 3 point
cases ($s=1,2$) introduced in Section 3 we abtain the following
deformed algebras,
\nl
(A) the full Virasoro algebra $\ \Cal V\ $ in all cases
of the cuspidal degeneration,
\nl
(B) the subalgebra $\ {\Cal W}^\a$, $\a=\sqrt{3e}\ $ of the
Virasoro algebra
in the case of the nodal degeneration
of the 2 point situation if $(e_1,0)$ does not become the singular point,
\nl
(C) the full vector field algebra  $\ {\Cal Z}^\a$, $\a=\sqrt{3e}\ $
 of 3 markings
in the case of the nodal degeneration either in the 2 point
situation
if $(e_1,0)$ becomes a singular point
or all 3 point situations with the exception of
the case considered in (D),
\nl
(D) a subalgebra of $\ {\Cal Z}^\a$, $\a=\i\sqrt{3e}\ $
(which could be given explictly) for the nodal degeneration in
the 3 point situation  with $e=e_1=e_2$ and $q=1/2+\t/4$.

\noindent All isomorphisms $\Phi $
induced by the degenerations
 respect the almost gradings of the
involved algebras.
\endproclaim
\np
{\bf 6. Further results}
\vskip 0.3cm
As already indicated, what have been done here for the vector field
algebra could have been done for the whole Lie algebra of
differential operators of degree less or equal one and their
Lie modules consisting of sections of $K^{\l}$ of arbitrary
integer weight $\l$. One obtains again by the well defined geometric
deformation certain subalgebras, resp\. modules of the
corresponding objects
on $\P^1$. Note however that for forms of positve weight
additional poles are introduced at the points lying above the
singularities.
Roughly speaking, one obtains by degenerations subalgebras
operating on bigger spaces.
This has to be considered if one uses this moduls to
construct semi-infinite wedge representations (i.e. $b-c$ systems),
see \cite{16,17}.

In this letter I only considered the algebras without central
extensions. As can easily seen the geometric cocycles depending
on the partition of the markings defining central
extensions \cite{16-18} can be calculated completely
outside the singularities. Moreover, they are calculated by
calculating residues (for example, at the point $\infty$),
hence in an algebraic manner.
They make sense under the described degeneration
process and in this way everything makes sense also for
the central extensions.

At least in principle the method is not restricted to the genus 1 case.
However, additional problems occur at higher  genus.
In the genus 1 case   I used that
every elliptic curve can be realized as the set of zeros
 of one polynomial
in the projective plane. The degenerations could be easily
studied at the defining polynomial. We obtained that there are
exactly two nonisomorphic possible degenerations.
In general, curves of higher genus can
 not be embedded in the complex plane.
One has to use projective spaces of higher dimensions
and needs more polynomial to define them. Under degeneration of the
curves
(if one considers only `stable' degenerations) one approaches
the boundary of the compactified moduli space of curves of
genus $g$. This boundary is of positive
dimension and consists of different  components, which can be
identified with the moduli space of lower genus curves
 (together with their
degenerations) \cite{19,p.78}.
Even without considering the markings
it is not at all clear in which `direction'
one should degenerate. A possibility could be `maximal degeneration'
in the sense that one obtains a (reducible) connected union of curves
either isomorphic to $\P^1$ or with desingularization $\P^1$.
By desingularization of the whole curve one obtains
a disjoint union of many copies of $\P^1$s.
Hence, one expects under the limit procedure for the
vector fields a subalgebra
of the direct sum of the associated vector field algebras on
different copies of $\P^1$s.
\np

\bye